
\input diagrams
%
%

\def\Z{{\rm Z}}
\def\N{{\rm N}}
\def\C{{\rm C}}

\def\R{{\rm R}}
\def\Q{{\rm Q}}
\def\SS{{\cal S}}
\def\Z2{{\rm Z}/_2}
\def\LM{LM\ }


\centerline{\bf Extended powers and Steenrod operations in algebraic geometry.}
\smallskip
\centerline{(Preliminary Draft, July 2007 version)}
\vskip .25 true cm
\centerline{Terrence Bisson \quad\& \quad Aristide Tsemo}
\centerline{bisson@canisius.edu \quad \  \quad tsemo58@yahoo.ca}
\vskip .25 true cm

\bigskip
\bigskip
\centerline{\bf Abstract.}
\smallskip

Steenrod operations have been defined by Voedvodsky in motivic cohomology in order to show the
Milnor and Bloch-Kato conjectures. These operations have also been 
constructed by Brosnan for Chow rings.
The purpose of this paper is to provide a setting for the construction of the Steenrod operations in
algebraic geometry, for generalized cohomology theories whose formal group law has order two. 
We adapt the methods used by 
Bisson-Joyal in studying Steenrod and Dyer-Lashof operations in unoriented cobordism and mod $2$
cohomology.

\centerline
{\bf Introduction.}
\smallskip

The mod $2$ cohomology ring $H^*(X;\Z2)$ of any space $X$ is naturally endowed with
operations; see Steenrod [1962]. The Steenrod square operations satisfy natural compatibility relations
such as the Adem relations which are complicated to state.
Bullet and McDonald [1982] (and Bisson [1977]) noticed that it is possible to formulate 
these relations in a convenient way, using formal power series.  The theory of $Q$-rings described 
in  Bisson, Joyal [1995a,b] [1997] incorporates this approach, and provides a setting
for Steenrod operations within an algebra of covering spaces,
interpreted as extended power functors in the category of topological spaces. 
Unoriented cobordism and the Thom realization functor transport the extended power functors
to give operations in $\Z2$-cohomology. In this setting the structure of 
$Q$-ring appears naturally, and then the proof of the Adem relations, and the rest of the theory, 
is straightforward.
 
Some of the ideas described by Bisson and Joyal were inspired by the paper of Quillen [1971], which
has also motivated Levine and Morel in their work on algebraic cobordism.  Let $k$ be a field and 
let $\SS$ be the category of quasi-projective schemes defined over $k$.
In the terminology of Levine and Morel, an {\it oriented cohomology theory}
on $\SS$ is a ring valued functor which satisfies various axioms. 
We will refer to these as \LM cohomology theories.
Over characteristic zero, Levine and Morel's algebraic cobordism 
is the universal example of this type of functor.
The purpose of this paper is to define extended power functors on $\SS$
as a setting for the construction of Steenrod operations.
In this we adapt methods from Bisson, Joyal [1995a].
As background, we note that the theory of mod 2 $Q$-rings is based 
on the fact that the mod $2$ cohomology 
of the topological classifying space of $\Z2$ is free on a formal variable $t$.
Classifying spaces for finite groups in algebraic geometry have
recently been defined by Morel and Voedvodsky and Totaro.
We will follow Totaro in working with certain affine schemes
built from respesentations as classifying space approximations.
In particular, for any finite group $G$ we will define
a sequence $B_nG$ of affine schemes
determined by the action of $G$ on the group algebra $k[G]$.

\medskip
Let $A$ be a \LM cohomology theory.
In order to define Steenrod-type operations, we need to make some additional assumptions 
on $A$:  

\smallskip\noindent
\item{$\bullet$} We assume that the formal group law $F_A(x,y)$ 
determined by $A$ satisfies $F_A(x,x)=0$, 
and similarly for the double covering formal group law.

\smallskip\noindent
\item{$\bullet$} We assume that $\lim_{n\to\infty}A(B_n\Z2)=A[[t]]$, 
the ring of formal power series over the coefficient ring for the theory $A$,
where $t$ is the characteristic class for double coverings.

\smallskip\noindent
\item{$\bullet$} We assume, for $G=\Sigma_4$, that 
any inner automorphism on $G$ induces the identity on $\lim_{n\to\infty}A(B_nG)$.

\smallskip\noindent
\item{$\bullet$} We assume the existence of a well-defined external 
extended power operation on $A$ that satisfies a few simple naturality conditions.

\smallskip\noindent
\item{$\bullet$} We assume that the resulting diagonal extended power operation on $A$ 
is additive for double coverings. 

We note that some of these assumption follow from axioms for oriented cohomology 
suggested by Panin and Smirnov [2000].

\bigskip

The use of methods of algebraic topology in algebraic geometry has a long history, including the work
of Grothendieck and his collaborators on defining a good framework for proving the Weil
conjectures. The approach of Grothendieck is wide-spread in algebraic geometry and has lead
to the proofs of many conjectures. Recently, Voedvodsky has shown the Milnor and Bloch-Kato conjectures
by using Steenrod operations in motivic cohomology. These operations have also 
been constructed in Chow rings by Brosnan. Our method is an attempt to situate their work in a
simple framework.

ACKNOWLEDGEMENT

The first author would like to thank Andr\'e Joyal for the gift of 
beautiful ideas in mathematics, and for our continuing collaborative work
on operations in algebraic topology.  This foray into the strange world of
algebraic geometry is entirely based on that work.

The first author also thanks Aristide Tsemo for initiating this enjoyable project,
and the participants at Canisius REU 2007 for keeping the summer lively.

The second author would like to thank Canisius College,  
in particular Terrence Bisson,
for warm hospitality. 
He wants to thank also his brother in law Maurice Neyou who enabled
him to meet Terrence Bisson.

\medskip
\noindent {\bf TABLE OF CONTENTS:}

1. Extended power functors in topology.

2. Background on algebraic geometry.

3. Extended power functors in algebraic geometry.

4. \LM cohomology theories.

5. Some axioms for extended power operations.

6. Some properties of extended power operations.

7. D-rings and mod 2 Steenrod operations.

\bigskip\noindent
{\bf 1. Extended power functors in topology.}

In the topological setting, a covering space is a continuous map
$p:T\to B$ which is locally trivial, with a finite number of sheets over each connected component.
Such a covering space can be used to define a functor 
from the category of topological spaces to itself.
This concept is developed and applied 
in two Comptes Rendus by Bisson and Joyal [1995a,b], and we will
closely follow that presentation here.
For any topological space $X$ we define
$$p(X)=\{(u,b)|b\in B, u:p^{-1}(b)\to X\}.$$
This construction is functorial for topological spaces and 
continous maps.
We will say that such a functor is an {\it extended power functor}.

The extended power functors could also be called polynomial functors.
Since the term ``extended power'' seems well established in topology, 
we have chosen that terminology here.

Here is another way of thinking about this extended power construction.
Let ${\bf n}=\{1,\ldots,n\}$, and let $\Sigma_n$ denote the group of bijections
of ${\bf n}$.  Suppose that all the fibers of the covering $p:T\to B$ have cardinality $n$;
for any $b\in B$, let ${\rm Frame}_b(p)$ denote the set of bijections 
from ${\bf n}$ to $p^{-1}(b)$.
The principle $\Sigma_n$-bundle $E(p)\to B$ associated to $p:T\to B$ has total space
$$E(p)=\{(b,f)|b\in B, f\in {\rm Frame}_b(p)\}$$
with $\Sigma_n$ acting freely on it.
Then $T=(E(p)\times {\bf n})/\Sigma_n$ and $p(X)=(E(p)\times X^n)/\Sigma_n$.
In this way we see that $p(X)$ is the total space of a bundle over $B$ with
fiber $X^{p^{-1}(b)}$ for $b\in B$.

Suppose that we think of isomorphism classes of topological spaces as forming a
``ring'' with disjoint union as $+$ and cartesian product as $\times$.
It is observed in Bisson, Joyal [1995a] that these
extended power functors are closed under these
operations of sum, product, and composition of functors 
from the category of topological spaces to itself.
In other words, given coverings $p$ and $q$ there exist
coverings $p+q$, $p\times q$, and $p\circ q$ such that we have, for $X$,
$$(p+q)(X)=p(X)+q(X)\quad{\rm and}\quad
(p\times q)(X)=p(X)\times q(X)\quad 
{\rm and}\quad (p\circ q)(X)=p(q(X)).$$

We can explicate these operations by using
the {\it derivative} $p'$ of a covering $p:T\to B$, which is defined
to be the covering with base space
$T$ such that the fiber of $p'$ over $t\in T$ is the set $p^{-1}(p(t)) - \{t\}$. 
We have rules like those for differential calculus:
$$(p+q)'=p'+q' \quad{\rm and}\quad (p\times q)'=p'\times q+p\times q'\quad{\rm and}
\quad(p\circ q)'=(p'\circ q)\times q'.$$
In fact, the covering $p$ can be expressed as $p'(1)\to p(1)$
(where $1$ denotes a singleton);
and the pull-back of $p$ along the natural map $p(X) \to p(1)$ can be identified
with a natural map $p'(X)\times X \to p(X)$, which is thus a covering.

Then the coverings $p+q$, $p\times q$, and $p\circ q$ are given by the natural maps
$$p+q:p'(1)+q'(1)\to p(1)\times q(1)\quad\quad
p\times q:p'(1)\times q(1)+p(1)\times q'(1)\to p(1)\times q(1)$$
$$\quad{\rm and}\quad p \circ q: p'(q(1))\times q'(1) \to p'(q(1))\times q(1) \to p(q(1)).$$

This shows that in the topological setting there is an algebra of coverings, 
wherein the sum, product, and composition
satisfy indentities appropriate to an algebra of polynomials.
Several applications of these concepts in algebraic topology are given 
in Bisson, Joyal [1995a,b].  One observation there is that 
if the base space of covering $p$ is a smooth manifold, then $X\mapsto p(X)$
is a functor from the smooth category to itself.

\smallskip
We want to transport the above concepts into the setting
of algebraic geometry, and show that they are relevant to the description of
Steenrod-type operations in algebraic geometry.

\bigskip\noindent
{\bf 2. Some background in algebraic geometry.}

It seems appropriate to work in the category of smooth quasi-projective schemes over a field.
We start by sketching some definitions and results from algebraic geometry. 

Let $k$ be a field.  
Each commutative $k$-algebra $R$ determines an {\it affine} $k$-scheme ${\rm Spec}(R)$.
The elements $x\in {\rm Spec}(R)$ correspond to the prime ideals of $R$; the
set of elements is given the Zariski topology and a distinguished sheaf of local rings over this topology.
Morphisms are defined so that the category of affine $k$-schemes 
is opposite to the category of commutative $k$-algebras.
The category of $k$-schemes, including notions of image and of open and closed 
subschemes, is modeled on this category.  Eisenbud and Harris [2000] give a very nice treatment.
They describe, for instance, how any $k$-scheme $X$ can be understood 
through its functor of points, which assigns to each $R$
the set $X(R)$ of scheme morphisms from ${\rm Spec}(R)$ to $X$.

As an example, let ${\rm A}^n$ denote {\it affine $n$-space} ${\rm Spec}(k[x_1,\ldots,x_n])$,
so that ${\rm A}^n(R)=\{(a_1,\ldots,a_n):a_i\in R\}$.  
Similarly, projective space ${\rm P^n}$ is the scheme whose functor of points
assigns to each $R$ the set of equivalence classes 
$[a_0,\ldots,a_n]$, where the $a_i\in R$ are not all zero and 
$[a_0,\ldots,a_n]\equiv[\lambda a_0,\ldots,\lambda a_n]$ for $\lambda\in k$ non zero.
A {\it projective} $k$-scheme is a closed subscheme of some projective space;
a {\it quasi-projective} $k$-scheme is any open subscheme in a projective $k$-scheme.
For instance, the affine schemes {\it of finite type} 
(opposite to the category of finitely generated commutative $k$-algebras)
are quasi-projective,
since each can be identified with an open subscheme in some projective scheme.

A morphism of $k$-schemes $f:X\rightarrow Y$ is a {\it closed imbedding} 
if and only if there exists a closed
subscheme $Y'$ of $Y$ such that $f$ factors by an isomorphism $X\rightarrow Y'$.
The morphism $f$ is a projective morphism, if and only 
if $f$ is the composition of a closed imbedding
$X\rightarrow {\rm P}^n_Y$ by the canonical surjection ${\rm P}^n_Y\rightarrow Y$
(here ${\rm P}^n_Y$ denotes the {\it relative projective $n$-space over $Y$}, 
which is often just ${\rm P}^n\times Y$).

Within the category of $k$-schemes, let $\SS$ denote the full
subcategory of smooth quasi-projective schemes.
The category $\SS$ has terminal object $1={\rm Spec}(k)$, 
and is closed under finite products and coproducts, denoted by $\times$ and $+$.
The existence of fiber products in $\SS$ depends on transversality.
The affine spaces and projective spaces
and their smooth closed subschemes are in $\SS$.
If $E\to B$ is a vector bundle in  $\SS$ then the associated projective bundle  
${\rm P}(E)\to B$ is in $\SS$.  Eisenbud and Harris [2000] give the relevant definitions.

Suppose that $E\to X$ is a vector bundle over $X$; then an $E$-torsor is a fiber bundle $S\to X$ 
together with a map $E\times_X S\to S$ over $X$ such that the associated map $E\times_X S\to S\times_X S$ 
is an isomorphism over $X$.
This gives a principal action of each fiber of $E$ on the corresponding fiber of $S$.
Since a principal action of a vector space on a set gives that set the structure of an affine space,
$S\to X$ may be called a {\it bundle of affine spaces} on $X$.

\bigskip

We end this section with some comments about the existence of 
categorical quotients by finite group actions.
If a finite group acts on a scheme $X$ then one would {\it like} to have a morphism
of schemes $X\to Y$ satisfying the universal property of categorical quotient.
This is impossible in general.
Even for a free action of a finite group, 
such a quotient of a scheme does not exist automatically in the category of schemes; 
see for instance the example of Hironaka [1962], or its description on page 15 in Knutson [1971].  

But for a free action of a finite group $G$ on a quasi-projective $k$-scheme $X$ there {\it does}
exist a $k$-scheme $X/G$ and a morphism of $k$-schemes $X\to X/G$ which is a categorical quotient.
This follows from the fact that any orbit of $G$ acting on quasi-projective $X$ is contained
in an affine open subscheme of $X$; see page 69 in Mumford [1970] (the argument there holds for
any field).  For context, see also the discussion around Proposition 1.8 in Expos\'e V
of Grothendieck [1971].

\bigskip

A {\it principal $G$-bundle} $E\to B$, for a finite group $G$, is just
a free action of $G$ on an affine scheme $E$ of finite type, with $B=E/G$.  
For example, let $\rho:G\times V\rightarrow V$ be a faithful representation of $G$
on a finite dimensional $k$-vector space $V$.  
The affine space $V$ has a closed subscheme $S$ on whose complement
$G$ acts freely; $S=\cup_{g\neq 1} V_g$, where $V_g$ is the linear subspace fixed by $g$.
This gives a free action of $G$ on $V-S$, a smooth affine scheme of finite type. 
Totaro [1999] shows that every principal $G$-bundle $E\to B$ (with $B$ an affine $k$-scheme)
is the pullback of one of these, for some $V$.  


\bigskip\noindent
{\bf 3. Extended power functors in algebraic geometry.}

The discussion of quotients from section 2 leads us to 
a convenient notion of covering spaces in algebraic geometry.

Any principal $\Sigma_n$ bundle $E\to B$ gives a morphism of affine schemes $p:T\to B$,
by taking $T=(E\times n)/\Sigma_n$, where $n=1+\cdots+1={\rm Spec}(k^n)$.
We will refer to such $p:T\to B$ as {\it geometric coverings with $n$ sheets};
a {\it geometric double covering} is just a geometric covering with $2$ sheets.

We can recover $E=E(p)$ from the geometric covering $p$ as in the topological setting.  
The scheme $T^n$ is affine, and the symmetric group $\Sigma_n$ acts naturally 
on $T^n$ (by permutating the indices of the $n$-tuples).
We consider the subscheme of all $n$-tuples $(t_1,\ldots ,t_n)$ in $T^n$ such that
$p(t_i)=p(t_j)$ for all $i,j$; this can be defined by repeatedly taking the fiber product
of $T\to B$ with affine schemes over $B$, and is thus afine.
The symmetric group $\Sigma_n$ acts naturally here, and 
$E(p)$ is defined to be the open subscheme with $t_i\neq t_j$ if $i\neq j$.
Then $E(p)$ is the total scheme of a principal $\Sigma_n$-bundle over $B$.

We say that $p:T\to B$ is a {\it smooth} geometric covering if and only if $T$ and $B$ are smooth;
by transversality of the fiber products defining $E(p)$,
this is equivalent to the condition that $E(p)$ be smooth.

Suppose that $p:T\to B$ is a smooth geometric covering with $n$ sheets.
For any smooth quasi-projective scheme $X$ we define $p(X)$ to be the quotient of 
$E(p)\times X^n$ by the diagonal action of $\Sigma_n$. 
The existence of this quotient is ensured by the result mentioned in section 2, 
since $E(p)\times X^n$ is a quasi-projective scheme.  Since $E(p)$ and $X$ are smooth,
so is $p(X)$.  We can use the following argument to show that  $p(X)$ is quasi-projective.
Totaro shows the existence of a $\Sigma_n$-equivariant 
closed embedding $E(p)\to V-S$ for some linear representation, as discussed above. 
But $X$ is an open subscheme of some projective $k$-scheme $\bar X$,
and the quotient of the $\Sigma_n$ action on the 
projective $k$-scheme ${\rm P}(V\oplus k)\times {\bar X}^n$ is a projective $k$-scheme.
The result follows.

Thus this construction defines
a functor $p:\SS\to \SS$ for each  covering in $\SS$.
Any functor $F:\SS\to \SS$ which is isomorphic, via a natural transformation, to
such a functor (for some $p$) is called an extended power functor. 

Given coverings $p$ and $q$ in $\SS$, we use the same formulations as for topological spaces to define
 coverings $p+q$, $p\times q$, and $p\circ q$.

\medskip
{\bf Proposition.}

{\it If $F$ and $G$ are extended power operations then the functors $F+G, F\times G, F\circ G$,
defined respectively by $(F+G)(X)=F(X)+G(X), (F\times G)(X)=F(X)\times G(X)$, and 
$(F\circ G)(X)=F(G(X))$, are extended power functors.}

The proof is the same as for topological spaces.


\bigskip
{\bf 4. \LM cohomology theories.}

Levine and Morel [2001a] introduced axioms for a notion of ``oriented cohomology theory'' 
in algebraic geometry.
These axioms are inspired by the method for developing complex cobordism theory which is 
presented in Quillen [1971].  
There Quillen suggests working with contravariant functors (from smooth manifolds to rings) 
which have covariant (or gysin) morphisms for proper smooth maps endowed with a complex orientation.
Multiplicative generalized cohomology theories 
which are oriented over complex cobordism provide examples for Quillen's discussion,
but Quillen does not require that his contravariant functors satisfy 
the full, usual axioms for a generalized cohomology theory.

Levine and Morel work on the category $\SS$.  
They assume the existence of gysin homomorphisms for a restricted category $\SS'$
of morphisms, those  which are projective morphisms of pure codimension.
A morphism $f:Y\rightarrow X$ in $\SS$
has {\it pure codimension} $d$ if we have ${\rm dim}_k(X,f(y))-{\rm dim}_k(Y,y)=d$ 
at every point $y$ in $Y$, 
where ${\rm dim}_k(Y,y)$ is the Krull dimension of $Y$ in a neighborhood of $y$.
Note that $\SS'$ contains the identity morphisms and is closed for composition,
do it does in fact form a subcategory of $\SS$ 
(with all smooth quasi-projective schemes as objects).

To allow for different conventions in handling dimensions, we will 
attach some fixed ``grade multiple'' $a$ to the theory; see LM1 below.
In our examples, $a$ is one or two.

Suppose that $a$ is a fixed integer, and $A$ is a contravaiant functor from $\SS$ to the 
category of graded commutative rings and grade-preserving ring homomorphisms;
a morphism of schemes $f$ gives a ring homomorphism $f^*$.
We will say that the pair $(A,a)$ is a  {\it \LM cohomology theory} 
if the functor $A$ and integer $a$ satisfy the following axioms
LM1-LM4 from Levine and Morel [2001a].

\medskip
{\bf LM 1:} $A$ is also a covariant functor from $\SS'$ to the category of graded abelian groups,
taking a morphism $f:Y\to X$ of pure codimension $d$ to a 
homomorpism $f_*:A(Y)\to A(X)$ which raises the grading by $ad$.
\medskip

From the contravariance we have a natural map $A(X)\times A(Y)\to A(X\times Y)$ for all $X$ and $Y$,
given by the multiplication in the ring $A(X\times Y)$.  

As a consequence of the covariance along $\SS'$, each projective morphism $f:Y\to X$ 
of codimension $d$ in $\SS'$ gives a class  ${\rm cl}(f)$ in $A^{ad}(X)$, defined by
${\rm cl}(f)=f_*(1)$ for $1\in A^0(Y)$.
Another consequence is the definition of an euler class for each vector bundle in $\SS$.
More precisely, if $\nu:E(\nu)\to X$ is a rank $n$ vector bundle in $\SS$, 
then the zero section $s:X\to E(\nu)$ is a smooth projective morphism of pure codimension $n$,
and we define $e(\nu)=s^*s_*(1)$ in $A^{an}(X)$, where $1$ is the identity in the ring $A(X)$.
It follows that the euler class satisfies $e(\nu_1\oplus \nu_2)=e(\nu_1)e(\nu_2)$ for vector bundles
$\nu_1$ and $\nu_2$ on $X$ in $\SS$.

For the next axiom, $\pi:P(\nu)\to X$ is the projective bundle 
of a rank $n$ vector bundle $\nu$ on $X$ in $\SS$, and $\gamma$ is the tautological 
line-bundle on $P(\nu)$.

\medskip
{\bf LM 2:}  $A(P(\nu))$ is a free $A(X)$-module 
with basis $1,e(\gamma),\ldots,e(\gamma)^{n-1}$ for every rank $n$ 
vector bundle $\nu:E(\nu)\to X$  in $\SS$  
\medskip

By methods of Grothendieck [1958], this axiom allows the definition of a 
complete family of characteristic classes for vector bundles, with results 
like those in Milnor, Stasheff [1974].
Also, since ${\rm P}(0)=\emptyset$ for the rank $0$ vector bundle, this axioms implies $A(\emptyset)=0$.

Let $f:Y\rightarrow X$ be a projective morphism in $\SS'$, and let $g:Z\rightarrow X$ be a morphism
in $\SS$ which is transverse to $g$, giving the scheme $Z':=Y\times_X Z$ in $\SS$ with projections
$f':Z'\rightarrow Z, g':Z'\rightarrow Y$, as in the following diagram.  
Then we say that $f$ and $g$ form a {\it transversal pullback} diagram.

\medskip
{\bf LM 3:} If $f$ in $\SS'$ and $g$ in $\SS$ form a transversal pullback diagram, 
then $f'_*\circ g'^*=g^*\circ f_*$
$${\rm if} \diagram
          Z' & \rTo^{f'} & Z  \cr
   \dTo_{g'} &           &  \dTo_{g}   \cr
          Y  & \rTo^f    & X
                                           \enddiagram \quad{\rm is\ transversal\ pullback,\ then}
\diagram
          A(Z') & \rTo^{f'_*} & A(Z)  \cr
   \uTo_{g'^*} &           &  \uTo_{g^*}   \cr
         A( Y)  & \rTo^{f _*}   & A(X)
                                           \enddiagram\quad {\rm commutes.}$$                   

In particular, this axiom computes $g^*({\rm cl}(f))={\rm cl}(f')$ for any projective morphism $f$ with
transversal $g$.  Also, we can deduce that $A(X+Y)=A(X)\oplus A(Y)$ (the coproduct of rings),
by applying this axiom to the transversal pullback diagrams given by $X\to X+Y$, $X\to X+Y$
and $X\to X+Y$, $Y\to X+Y$.

The next axiom is a partial ``homotopy'' axiom, among other consequences.  

\medskip
{\bf LM 4:} If $\nu:E\rightarrow X$ is a vector bundle over $X$ in $\SS$, 
then $\nu^*:A(X)\rightarrow A(E)$ is an isomorphism, and the same for any 
bundle of affine spaces on $X$ in $\SS$.

\medskip
A formal group law defined over a commutative ring $R$ is a formal power series $F(x,y)\in R[[x,y]]$
 which satisfies identities corresponding to associativity and unit and inverses.
Levine and Morel [2001a] explain how to deduce from their axioms
 the existence of a formal group law $F(x,y)$ with coefficients in the ring $A=A(1)$, such that
$e(\gamma_1\otimes\gamma_2)=F(e(\gamma_1),e(\gamma_2))$ for all line bundles 
$\gamma_1,\gamma_2$ on $X$ in $\SS$.

\bigskip
In their monograph [2007], Levine and Morel present the following examples 
(and others) which satisfy their axioms:
\medskip
1. The functor which sends a quasi-projective scheme $X$ defined over the field $k$ to the 
Chow ring $CH^*(X)$ is an oriented cohomology theory.
\medskip
2. Let $\ell$ be a prime number distinct from the characteristic of $k$; the functor which sends
$X$ to the sum of etale groups $\oplus_nH^{2n}_{et}(X,{\Q}_l(n))$ is an oriented cohomology theory.
\medskip
3. Let $K^0(X)$ be the Grothendieck group of vector bundles on the scheme $X$; 
the functor which sends $X$ to the
ring of Laurent series $K^0(X)[\beta,\beta^{-1}]$ is an oriented cohomology theory.
\medskip
4. Let $k$ be a number field and $\sigma:k\rightarrow {\C}$ be a complex embedding. For each
quasi-projective scheme $X$, we denote by $X_{\sigma}({\C})$ the quasi-projective variety of
complex points defined by $\sigma$. Let $MU$ be the complex cobordism spectrum, the functor
$X\rightarrow MU(X_{\sigma}({\C}))$ is an oriented cohomology theory.

\bigskip
When the oriented cohomology theory $A$ is very closely related to ordinary cohomology,
the formal group law may be the additive formal group law defined by $F(x,y)=x+y$, but in
general $F$ is more complicated. For instance example 3 has the
formal group law $F(x,y)=x+y-\beta x y$.

\bigskip
{\bf 5. Some axioms for extended power operations.}

Let $p:T\to B$ be an $n$-sheet geometric covering in $\SS$.
Consider an arbitrary \LM cohomology theory $(A,a)$ on $\SS$.
Each codimension $d$ projective morphism $f:Y\to X$ in $\SS'$ 
represents a cohomology class ${\rm cl}(f)\in A^{ad}(X)$.
Then the extended power functor $p$ gives $p(f):p(Y)\to p(X)$,
which represents a cohomology class ${\rm cl}(p(f))\in A^{nad}(p(X))$
(an argument similar to those in section 4 shows that $p(f)$ is a projective morphism
of codimension $nd$).   This {\it suggests} that the extended power functor $p$ may give
an ``external'' cohomology operation from $A^{ad}(X)$ to $A^{nad}(p(X))$,
and we can use geometric calculations in $\SS$ to guess at properties that such a
cohomology operation would have.

Unfortunately, we do not know that every class in $A(X)$ is represented by an $f\in \SS'$, 
and we have not shown that ${\rm cl}(f)={\rm cl}(f')$ implies ${\rm cl}(p(f))={\rm cl}(p(f'))$.

It seems reasonable at this stage to introduce additional axioms that a \LM theory
should satisfy, if it is to be equipped with extended power operations underlying a notion
of Steenrod operations.  That is the purpose of this section.  Then in the next
section we define Steenrod operations in such a cohomology theory, and develop 
their basic properties.  We limit ourselves to the case of $\Z2$
Steenrod operations in this article.

\bigskip
Let $(A,a)$ be a \LM cohomology theory in $\SS$, the category of quasi-projective schemes 
over a field $k$.  We assume for the rest of this paper that the 
following {\it Extended Power} axioms
(EP1, EP2, EP3, EP4, and EP5) are satisfied.

\medskip\noindent
{\bf EP1:} For every  cover $p:T\rightarrow B$ with $n$ sheets,
there exists a multiplicative map $p^{ext}:A^d(X)\rightarrow  A^{nd}(p(X))$ ({\it not} 
assumed to be additive in general) which:

a) is natural with respect to $X$: $p(f)^*\circ p^{ext}=p^{ext}\circ f^*$,

b) agrees with the $n^{th}$-power map when $p$ is the trivial geometric covering with $n$ sheets,

c) is natural with respect to $p$, and

d) commutes with euler classes: $p^{ext}(e(\nu))=e(p(\nu))$,

\medskip\noindent
Let us make this precise.

\medskip
\item{}For a), any $f:X\to Y$ in $\SS$ gives $p(f):p(X)\to p(Y)$ in $\SS$, 
and we require that $p(f)^*\circ p^{ext}=p^{ext}\circ f^*$.

\medskip
\item{}For b), we require that $p^{ext}(a)=a^n$ for $p:n\rightarrow 1$,

\medskip\noindent
For any $F:B'\rightarrow B$ and any $n$-sheeted geometric covering $p:T\rightarrow B$ in $\SS$, 
the pullback of $p$ along $F$ is a geometric covering $q:T'\to B'$ with $n$ sheets in $\SS$,
and we have a natural transformation $F(X):q(X)\to p(X)$ for each $X$.

\item{}For c) we require that
$${\rm if} \diagram
          T' & \rTo & T  \cr
   \dTo_{q} &           &  \dTo_{p}   \cr
          B'  & \rTo^F    & B
                                           \enddiagram \quad\quad{\rm is\ a\ pullback,\ then}
\diagram
          A^d(X) & \rTo^{q^{ext}} & A^{nd}(q(X))  \cr
               &  \SE_{p^{ext}} &  \dTo_{F(X)^*}   \cr
               &        & A^{nd}(p(X)) 
             \enddiagram\quad\quad {\rm commutes.}$$           

\item{}For d), let $\nu:V\rightarrow X$ be a rank $v$ vector bundle in $\SS$,
with Euler class $e(\nu)\in A^{an}(X)$.
If $p:T\rightarrow B$ is a geometric   covering 
with $n$ sheets then we have the rank $nv$ vector bundle $p(\nu):p(V)\rightarrow p(X)$.
We require that $p^{ext}(e(\nu))=e(p(\nu))$.

\bigskip
For any $n$ sheeted geometric covering $p:T\to B$ in $\SS$ 
we have a diagonal map $\Delta:p(1)\times X\to p(X)$ in $\SS$
(from the $\Sigma_n$ equivariant map $E(p)\times X\to E(p)\times X^n$).
This gives a diagonal pullback $\Delta^*:A(p(X))\to A(p(1)\times X)$, natural in $X$.
Assuming EP1, we may define $p^{\Delta}=\Delta^*\circ p^{ext}$.
Since our goal is $\Z2$ Steenrod operations, we make the following assumption.

\medskip
{\bf EP2:}
We assume that the map $p^{\Delta}:A^d(X)\rightarrow A^{2d}(p(1)\times X)$ is an additive homomorphism 
whenever $p$ is a geometric covering with two sheets.
\bigskip

From assumption EP1, the trivial geometric covering $p:2\to 1$ gives the extended power operation
$$p^\Delta:A(X)\to A(X^2)\to A(1\times X)=A(X)\quad a\mapsto a^2.$$
So in particular, EP2 implies that squaring is additive on $A(X)$.
This implies that $A(X)$ is always a ring of characteristic $2$.
In fact, we want to make a much stronger assumption.

Let $(A,a)$ be a \LM cohomology theory on $\SS$.  Let $F_A(x,y)$ be the
formal group law determined by $A$ (see the discussion after axiom LM4).
We say that a formal group law $F(x,y)$ in $A[[x,y]]$ 
has {\it order two} if $F(x,x)=0$ in $A[[x]]$.  
We say  that such a formal group law $F(x,y)$ in $A[[x,y]]$ 
is {\it compatible with $(A,a)$} if there exists $\tilde F(x,y)$, 
a formal group law of order two in $A[[x,y]]$, 
with $F(x^a,y^a)=(\tilde F(x,y))^a$ (this condition is vacuous if $a=1$).
We make the following assumption. 

\medskip
{\bf EP3:} We assume that the formal group law $F_A(x,y)$ determined by $A$ 
has order two, and also that $F_A$ is compatible with $(A,a)$.
\medskip

The condition $F(x^a,y^a)=(\tilde F(x,y))^a$ says that the power series $h(x)=x^a$ in $A[[x]]$
is a morphism of formal group laws $h:\tilde F\to F$; see Quillen [1971] and section 7 here.

\bigskip

For any finite group $G$, let $\overline{k[G]}$ denote
the reduced regular representation of $G$,
the kernel of the augmentation ring homomorphism
$\epsilon: kG\to k$ (with $\epsilon(g)=1$ for each $g\in G$). 
For $n>0$ let $n\ \overline{k[G]}$ denote the direct sum of 
$n$ copies of this representation,
and the corresponding affine space with its $G$ action.
Consider the open subscheme $(n\ \overline{k[G]}-S)$ where $G$ acts freely, and 
let $B_nG=(n\ \overline{k[G]}-S)/G$, the base of the
corresponding principal $G$ bundle.

For any injective group homomorphism $\phi:H\to G$ of finite groups,
we have an $H$-equivariant linear map $\phi: \overline{k[H]}\to \overline{k[G]}$;
since $\phi$ carries the free part $n\ \overline{k[H]}$ into the free part of
$n\ \overline{k[G]}$ (as we see by decomposing $k[G]$ along the cosets of $H$ in $G$), 
$\phi$ induces a morphism $B_n\phi:B_nH\to B_nG$ of affine $k$-schemes.

We use the following assumption  in our proof of the Adem relations. 

\medskip
{\bf EP4:} We assume that if $\phi:G\to G$ is an inner automorphism 
then $B_n\phi:B_nG\to B_nG$
gives the identity map on $A(B_nG)\to A(B_nG)$. 
\bigskip

A double covering $p:T\to B$ in $\SS$ determines a line bundle $\gamma(p)$ on $B$;
the total space of $\gamma(p)$ can be described as $(E(p)\times {\rm A}^1)/(\Z2)$,
where $E(p)$ is the principal $\Z2$ bundle for $p$, and $\Z2$ acts antipodally on ${\rm A}^1$.
A {\it characteristic class for double coverings} assigns a class $t(p)\in A^1(B)$
to each geometric double covering $p:T\to B$ in $\SS$, so that the assigment is
natural in $p$ and satisfies $e(\gamma(p))=t(p)^a$.  Note that if $a=1$, then
$t(p)=e(\gamma(p))$ determines such a characteristic class for double coverings.
In general, the existence of a characteristic class for double coverings is linked to the behavior
of the cohomology functor $A$ on a classifying space for principle $\Z2$ bundles.

For the rest of the paper 
we assume that our field $k$ is not of characteristic $2$.
Then every finite dimensional vector space $V$
provides a faithful representation of the group $\Z2$ acting as the antipode map $v\mapsto -v$;
this gives us the geometric  covering $p_V:(V-0)\to (V-0)/(\Z2)$ in $\SS$.
Let $p_n$ be the geometric  covering from the vector space $k^n$, 
and let $B_n\Z2$ denote the base $p_n(1)$
of this covering. This agrees with our above notation $B_nG$ since, as $\Z2$ representations,
$k$ with its antipode action is  isomorphic to $\overline{k[\Z2]}$.

We make the following assumption.

\medskip
{\bf EP5:} We assume the existence of a characteristic class for double coverings,
such that $A(B_n\Z2)=A[t]/t^{an}$, and
$A(B_n\Z2\times X)=A(X)[t]/t^{an}$ naturally in $X$, where $t=t(p_n)\in A^1(B_n\Z2)$.
We also assume that $p_n^\Delta:A^1(B)\to A^2(B_n\Z2\times B)$ satisfies
$p_n^\Delta(u)=u{\tilde F}_A(u,t)$ whenever $u=t(p)$ is the characteristic class 
of a geometric double covering $p:T\to B$.

\medskip
If we assume that our 
cohomology theory satisfies the long exact sequence (excision axiom) 
from Panin and Smirnov [2000],
then we can derive a gysin sequence for vector bundles, and use this to compute
$A(B_n\Z2\times X)=A(X)[t]/t^{an}$, at least for $a=1$ and $a=2$.
The proof uses the fact that $B_n\Z2$ is isomorphic to the complement of the zero section
of a line bundle $\gamma\otimes\gamma$ over the projective space ${\rm P}^{n-1}$.  

\bigskip
Here are some natural examples where the extended power assumptions are satisfied.
Suppose that the field $k=\R$ is the
real numbers.  For any $X$ in ${\cal S}_{\R}$, the set $X(\R)$ of real-valued points in $X$
is a smooth manifold.  Let $N(X)={\cal N}^*(X(\R))$, the unoriented cobordism ring
of the smooth manifold $X(\R)$.  Let $H(X)=H^*(X(\R);\Z2)$, the mod $2$ cohomology ring
of the smooth manifold $X(\R)$.  Then $(N,1)$ and $(H,1)$ are \LM cohomology theories
which satisfy EP1, EP2, EP3, EP4, and EP5.
Note that $F_n$ is the additive formal group law, $F_H(x,y)=x+y$;
 and $F_N$ is the usual formal group law for unoriented cobordism, 
which is the universal formal group law of order $2$ (see Quillen [1971]).

Suppose instead that the field $k=\C$ is the
complex numbers.  For any $X$ in ${\cal S}_{\C}$, the set $X(\C)$ of complex-valued points in $X$
is a smooth even dimensional manifold.  
Let $NC(X)={\cal N}^*(X(\C))$ and let $HC(X)=H^*(X(\C);\Z2)$.
Then $(NC,2)$ and $(HC,2)$ are \LM cohomology theories
which satisfy EP1, EP2, EP3, EP4, and EP5.
For $(NC,2)$, we note that
the formal group law $F_{NC}(x,y)$ is determined by the tensor product of complex line bundles 
rather than the tensor product of real line bundles.  
In fact, $F_{NC}(x^2,y^2)=(F_N(x,y))^2$.

\bigskip\noindent
{\bf 6. Some properties of extended power operations.}

Recall that $k$ is a field of characteristic
different from $2$, and that $\SS$
is the category of quasi-projective schemes over $k$.
We assume for the rest of the paper
that $(A,a)$ is a \LM cohomology theory with
extended power operations satisfying the assumptions from the previous section.

Our first goal is to define a ``total operation'' for $A$.
Consider the geometric double coverings $p_n:T_n\Z2\to B_n\Z2$ determined
by the free action of $\Z2$ on $k^n-0$ (as discussed in connection with
assumption EP4). 
Note that $p_n$ is the pull-back of $p_{n+1}$ along the map $i_n:B_n\Z2\to B_{n+1}\Z2$
induced by the inclusion of $k^n$ into $k^{n+1}$ by the first $n$ coordinates.
This gives natural transformations such that the following diagram commutes:
$$\diagram
          p_n(X) & \rTo_{i_n} & p_{n+1}(X)  \cr
   \uTo_{\Delta} &           &  \uTo_{\Delta}   \cr
          p_n(1)\times X  & \rTo^{i_n\times id_X}    & p_{n+1}(1)\times X
\enddiagram$$
By assumption EP1 we have $i_n^*\circ p_{n+1}^{ext}=p_n^{ext}$. Composing with
$\Delta^*$ gives $i_n^*\circ p_{n+1}^{\Delta}=p_n^{\Delta}$.
By assumption EP5 we identify $i_n^*:A(p_{n+1}(1)\times X)\to A(p_{n}(1)\times X)$
with the truncation $A(X)[t]/t^{a(n+1)}\to A(X)[t]/t^{an}$.
But we have an isomorphism 
$$\lim_{n\rightarrow \infty}A(B_n\Z2\times X)\rightarrow A(X)[[t]].$$
It follows that the limit of the operations $p_n^{\Delta}$ defines a total operation
$$P_t:A(X)\rightarrow A(X)[[t]].$$
Assumption EP2 implies that this is a ring homomorphism,
and we have the following.

\medskip
{\bf Proposition.}

{\it The functor $A$ on $\SS$ is equipped with a natural ring homomorphism
$P_t:A(X)\to A(X)[[t]]$.}

\bigskip

We will refer to $P_t$ as the {\it total double covering operation} for $A$.

\medskip
{\bf Proposition.}

{\it If $t$ is the characteristic class of a geometric double covering in $\SS$, then
$P_u(t)=t{\tilde F}_A(t,u)$.

If $e$ is the euler class of a line bundle in $\SS$
 then $P_u(e)=eF_A(e,u^a)$.}
 
 {\bf Proof.}

The first statement comes by taking the limit of the formula $p_n^\Delta(t)=t\tilde F(t,u)$ from
assumption EP5.  For the second statement,
let $\nu:V\rightarrow X$ be a rank $v$ vector bundle in $\SS$,
and let $e(\nu)\in A^n(X)$ denote the Euler class of $\nu$. 
If $p:T\rightarrow B$ is a geometric   covering 
with $n$ sheets then we can define the vector bundle $p(\nu):p(V)\rightarrow p(X)$; it has rank $nv$.
Consider $\Delta^*p(\nu)$, the pullback of $p(\nu)$ along the diagonal map 
$\Delta:p(1)\times X\rightarrow p(X)$. We have $p_{\Delta}(e(\nu))=e(\Delta^*(p(\nu)))$.
This follows from EP1.
We get a rank $n$ vector bundle $\rho$ on  $p(1)$ by applying $p$ to the one-dimensional $k$ vector space
(viewed as a trivial vector bundle on $1$).
We have $\Delta^*p(\nu)=\rho\otimes\nu$. Indeed these bundles have the same rank, and there 
exists a canonical morphism of bundles $\rho\otimes\nu\rightarrow \Delta^*p(V)$. We obtain that
$e(\Delta^*p(\nu))=e(\rho\otimes \nu)$. 
Suppose now that $\nu$ is a line bundle and that $p$ has two sheets; then we can decompose
$\rho=k\oplus \gamma$ for a line bundle $\gamma$ on $p(1)$. 
We deduce that $e(\rho\otimes\nu)=e(\nu\oplus(\gamma\otimes\nu))=e(\nu)e(\gamma\otimes\nu)$. 
Since $e(\gamma)=t^a$ in $A^a(p(1))$, this implies the result.  DONE.

\bigskip
Next we define the total operation for any finite group $G$.  

Totaro [1999] observes that any faithful representation of $G$ 
in a finite dimensional vector space $V$ over $k$
determines an geometric   covering in $\SS$, and then
$V$ gives a sequence of faithful representations $nV$ and a sequence
of  coverings.  In connection with our discussion of EP4, we defined
$p_{nG}:T_nG\to B_nG$ to be the geometric covering determined 
by the reduced regular representation $\overline{k[G]}$ of $G$.  
The number of sheets of $p_{nG}$ is the cardinality of $G$. 
We use geometric coverings $p_{nG}$ to define a total operation 
$$P_{nG}:A(X)\rightarrow \lim_{n\rightarrow \infty}A(B_nG\times X)$$
as  limit of the extended power operations $p_{nG}^\Delta:A(X)\to A(B_nG\times X)$.

When the order of $G$ is prime to the characteristic of $k$, then
Morel and Voevodsky [2001] describe a classifying object $B_{et}G$ for 
principal $G$ bundles in algebraic geometry. 
We may view the base spaces $B_nG$ of our geometric coverings $p_{nG}$ as 
affine $k$-scheme approximations to $B_{et}G$
(even when the order of $G$ is not prime to the characteristic).

\medskip
Recall that if $\phi:H\to G$ is an injective homomorphism of finite groups,
then we obtain a map $B_n\phi:B_n{H}\rightarrow B_n{G}$ (see the discussion of EP4), and a map
$$B\phi^*:
\lim_{n\rightarrow \infty}A(B_n{G}\times X)\rightarrow 
\lim_{n\rightarrow \infty}A(B_n{H}\times X)$$ 
such that $P_{H}=B\phi^*\circ P_{G}$.

\medskip
When $G$ is a subgroup of $\Sigma_r$, we say that 
a geometric covering $T\to B$ with $r$ sheets admits a {\it reduction of structure group} 
to $G$ if and only if
$p$ can be expressed in the form $T=(E\times r)/G$ for some principal $G$ bundle $E\to B$.

We need to analyze the following example.
Let $K$ be the subgroup of $\Sigma_4$ generated by the permutations $\tau_1=(12)(34)$ and 
$\tau_2=(14)(23)$ (the ``Klein four subgroup'' of $\Sigma_4$). 
The group $K$ is isomorphic to $\Z2\times \Z2$. The definition of the 
classifying space shows that $B_nK=B_n\Z2\times B_n\Z2$. This implies (by assumption EP5) that 
$$A(B_nK\times X)=
A(B_n\Z2\times B_n\Z2\times X)=A(B_n\Z2\times X)[u]/(u^{an})=A(X)[u,v]/(u^{an},v^{an}).$$
We deduce that the total operation $P_K$ is of the form
$$P_K=P_{u,v}:A(X)\rightarrow A(X)[[u,v]].$$

\medskip
{\bf Proposition.}

{\it  The map $P_{u,v}$ is symmetric, that is $P_{u,v}=P_{v,u}$.}

{\bf Proof.}

Consider the permutation $\tau=(1234)$ of $\Sigma_4$. For $\tau_1=(12)(34)$ and $\tau_2=(14)(23)$ 
as above we have 
$\tau\circ\tau_1\circ\tau^{-1}=\tau_2$. The embeddings $j_i:\Z2\rightarrow \Sigma_4,\ i=1,2$
induced by $\tau_1$ and $\tau_2$ give an embedding $j:K\rightarrow \Sigma_4$.
Since $\tau$ is an inner automorphism of $\Sigma_4$, assumption EP4 says that
the automorphism induced by $\tau$ on $B_n{\Sigma_4}$ is the identity.
We analyze the two maps $Bj_1,Bj_2:B_nK\to B_n\Sigma_4$ as above.
They are related by $\tau$.
This shows that $\tau$ induces on $B_nK$ an
automorphism which exchanges $u$ and $v$ in $A(X)[u,v]/(u^{an},v^{an})$. Since 
$P_K=Bj^*\circ P_{\Sigma_4}$, we deduce that $P_{u,v}$ is symmetric.  DONE.

\bigskip

We have defined $P_{u,v}:A(X)\to A(X)[[u,v]]$.  But this can also be interpreted
as an iteration of total double covering operations, as follows
(recall the definition of $\tilde F_A$ in connection with EP5).

\medskip
{\bf Proposition.}

{\it $P_{u,v}=P_u\circ P_v$ if we put $P_u(v)=v \tilde F_A(v,u)$.}

{\bf Proof.}

Let $K=\Z2\times\Z2$, as above.
A $K$-covering is a principal $K$-bundle $p:E\rightarrow E/K=B$, determined by a pair of 
fixed point free involutions $\tau_1,\tau_2:E\rightarrow E$ such that 
$\tau_2\circ\tau_1=\tau_1\circ \tau_2$. These determine a pair of double coverings 
$p_i:E_i\rightarrow B i=1,2$ where $E_i=E/\tau_i$. Let $p_1\otimes p_2$ be the $K$-covering
$E_1\times E_2\rightarrow B\times B$. The bundle $p$ is isomorphic to the pullback of 
$p_1\otimes p_2$ along the diagonal map $B\rightarrow B\times B$, and $p_1\otimes p_2$ is
isomorphic to the pullback of $p_1\circ p_2$ along the generalized diagonal map 
$p_1(1)\times p_2(1)\rightarrow p_1(p_2(1))$. For any $K$-covering $p$, let $u$ and $v$
in $A(p(1))$ be the characteristic classes 
of the two double coverings $p_1$ and $p_2$ associated to $p$.
The pullback arrow $\phi:p\rightarrow p_1\otimes p_2$ determines a ring homomorphism 
$\phi^*:A(p_2(1)\times p_1(1)\times X)\rightarrow A(p(1)\times X)$. We have 
$p=\phi^*\circ p_2\circ p_1$. But $p_1:A(X)\rightarrow A(p_1(1)\times X)$ comes from
$P_u:A(X)\rightarrow A(X)[[u]]$, and $p_2:A(p_1(1)\times X)\rightarrow 
A(p_2(1)\times p_1(1)\times X)$ comes from $P_v: A(p_1(1)\times X)\rightarrow
A(p_1(1)\times X)[[v]]$, while $p:A(X)\rightarrow A(X)[[u,v]]$ comes from
$P_{u,v}:A(X)\rightarrow A(X)[[u,v]]$. DONE.
 
\bigskip\noindent
{\bf 7. D-rings and mod 2 Steenrod operations.}

The following is taken from Bisson, Joyal [1995a], where the motivation was the study of
unoriented cobordism and bordism operations in topology.

 Let $R$ be a commutative ring and let $F(x,y)\in R[[x,y]]$ be a formal group law of order two
(note that this implies that $R$ is a $\Z2$ algebra). 
According to Lubin [1967], there exists a unique formal group 
law $F_t$ defined over $R[[t]]$ such that $h_t(x)=xF(x,t)$ is a morphism $h_t:F\rightarrow F_t$.
See Bisson, Joyal [1995a] for more discussion of this notion.
\bigskip
{\bf Definition.}

A $D$-ring $R$ is a commutative $\Z2$ ring endowed with a formal group law $F$ of order two

\noindent and with a total operation $D_u:R\rightarrow R[[u]]$ which satisfies the following conditions: 

\noindent D1: $D_0(a)=a^2$ for every $a\in R$;

\noindent D2: $D_u(F)=F_u$;

\noindent D3: $D_u\circ D_v$ is symmetric in $u$ and $v$, where
$D_u$ is extended to $R[[v]]$ by setting $D_u(v)=vF(u,v)$.

If $R^*$ is also a graded ring with $d_i(x)\in R^{2q-i}$ for $x\in R^q$,
then we say that $R$ is a {\it graded $D$-ring}.

\medskip
{\bf Theorem.}

{\it If $A$ is a \LM cohomology theory which satifies the extended power assumptions 
EP1, EP2, EP3, EP4, and EP5 then $P_u:A(X)\to A(X)[[u]]$, together with 
the formal group law $F=\tilde F_A$,
defines on $A(X)$ the structure of a graded $D$-ring,
for every $X$ in $\SS$.}

{\bf Proof.}

The fact that $P_0(u)=u^2$ is a direct consequence of assumption EP1. 
The formal group law $F=\tilde F_A$ has order 2 by 
assumption EP3.  Finally, we show that $P_u(F)=F_u$. Let $\gamma_1,\gamma_2$ be line bundles.  
Let $h_u(x)=xF(x,u)$.  We have shown, in the preceding section, that $P_u(e_i)=h_u(e_i),\ i=1,2$. 
Let $e$ be the euler class of the tensor product of
bundles $\gamma_1\otimes \gamma_2$, we have $F(e_1,e_2)=e$. Write $F(x,y)=\sum_{i,j} a_{ij}x^iy^j$;
the naturality of $P_u$ implies  that 
$$
\sum_{i,j}P_u(a_{ij})h_u(e_1)^ih_u(e_2)^j=h_u(F(e_1,e_2)).
$$
Since $F_u$ is the formal group law defined by $h_u(F(e_1,e_2))=F_u(h_u(e_1),h_u(e_2))$, we 
deduce that $P_u$ takes the coefficients of $F$ to the coefficients of $F_u$.  DONE.

\bigskip

When the formal group law is additive the notion of $D$-ring becomes the notion of $Q$-ring, as
defined in Bisson, Joyal [2001].
We have the following:

\medskip
{\bf Proposition.}

{\it If $A$ is a \LM cohomology theory which satifies the extended power assumptions 
EP1, EP2, EP3, EP4, and EP5 and $F_A(x,y)=x+y$ then $P_u$
defines on $A(X)$ the structure of a graded $Q$-ring
for every $X$ in $\SS$.}

\bigskip
The $Q$-ring structure $Q_u:R\to R[[u]]$ on the graded ring $R=A(X)$ is
equivalent to a sequence of individual additive operations $q_i:R^n\to R^{2n-i}$
determined by $Q_u(a)=\sum_{n\in{\N}}q_i(a)u^i$, when these operations
satisfy the $Q$-ring versions of the Cartan formula and the Adem relations;
see Bisson, Joyal [2001].  These individual operations capture exactly the structure
of $R$ as an unstable algebra over the Steenrod operations
${Sq}^j(x_n)=q_{n-j}(x_n)$ for $x_n\in R^n$.

This is the sense in which our extended power axioms ensure the existence of
Steenrod operations in an \LM cohomology theory in algebraic geometry.
\bigskip

This paper was typeset using Paul Taylor's \TeX  macros for diagrams.

\bigskip
\centerline{\bf REFERENCES}
\medskip
\noindent[1977] T. Bisson; Divided sequences and bialgebras of homology operations, Ph D. Thesis,
Duke University.
\smallskip
\noindent[1995a] T. Bisson, A. Joyal; The Dyer-Lashof algebra in bordism. C.R. Math. Rep. Acad. Sci. Canada
Vol. 17 n. 4, pp. 135-140.
\smallskip
\noindent[1995b] T. Bisson, A. Joyal; Nishida relations in bordism and homology. C.R. Math. Rep. Acad. Sci.
Canada. Vol 17 n. 4, pp. 141-145.
\smallskip
\noindent[1997] T. Bisson, A. Joyal; $Q$-rings and the homology of the symmetric groups. 
Contemp. Math., no. 202, Amer. Math. Soc., Providence, RI.
\smallskip
\noindent[2003] P. Brosnan, Steenrod operations in Chow theory.  Trans. Amer. Math. Soc.  355  
  no. 5, pp. 1869-1903. 
\smallskip
\noindent[1982] S. Bullet, G. McDonald; On the Adem relations. Topology 21 pp. 329-332.
\smallskip
\noindent[2000] D. Eisenbud, J. Harris; The geometry of schemes. Graduate texts in Mathematics
no. 197.
\smallskip
\noindent[1958] A. Grothendieck; La thŽorie des classes de Chern.   Bull. Soc. Math. France  86  
pp. 137-154.
 \smallskip
\noindent[1971] A. Grothendieck; Lecture Notes in Mathematics, no. 224. 
\smallskip
\noindent[1962] H. Hironaka; An example of a non-Khalerian deformation. 
Annals of Math. 75 pp. 190-208. 
\smallskip
\noindent[1971] D. Knutson; Algebraic spaces. Lecture Notes in Mathematics, no. 203. 
\smallskip
\noindent[2001a] M. Levine, F. Fabien; Cobordisme algebrique I.  C. R. Acad. Sci. Paris SŽr. I Math.  332  
 no. 8, pp. 723-728.
 \smallskip
\noindent[2001b] M. Levine,  F. Morel;  Cobordisme algŽbrique. II.  [Algebraic cobordism. II]  
C. R. Acad. Sci. Paris SŽr. I Math.  332    no. 9, pp. 815-820.
 \smallskip
\noindent[2007] M. Levine, F. Morel; Algebraic cobordism. Monographs in Mathematics. 
 \smallskip
\noindent[1967] J. Lubin; Finite subgroups and isogenies of one-parameter formal Lie groups,
 Ann. of Math. 85 pp. 296-302.
\smallskip
\noindent[1974] J. Milnor, J. Stasheff; Characteristic classes, Princeton University Press.
\smallskip
\noindent[2001] F. Morel, V. Voevodsky; $A^1$-homotopy of schemes.   Inst. Hautes ƒtudes
 Sci. Publ. Math.  No. 90   pp. 45-143.
\smallskip
\noindent[1970] D. Mumford; Abelian varieties. Tata Institute of Fundamental research. no. 5
\smallskip
\noindent[2000] I. Panin, A. Smirnov; Push-forward in oriented cohomology theories of algebraic
varieties. K-theory preprint archives.
\smallskip
\noindent[1971] D. Quillen; Elementary proofs of some results of cobordism theory using Steenrod operations,
Adv. in Math. 7  pp. 29-56. 
\smallskip
\noindent[1962] N. Steenrod, D. Epstein; Cohomology operations, Ann. Math. Studies, 
vol 50, Princeton University
Press. 
\smallskip
\noindent[1999] B. Totaro; The Chow ring of a classifying space. 
Proc. Sympos. Pure Math no. 67, Amer. Math. Soc.
pp. 249-281.
\smallskip
\noindent[2003] V. Voevodsky;  Reduced power operations in motivic cohomology.
  Publ. Math. Inst. Hautes ƒtudes Sci.  no. 98  pp. 1-57.

\bigskip\noindent
Terrence P. Bisson, Canisius College, Buffalo, NY 14216 USA

\bigskip\noindent
Aristide Tsemo, 44 St. Dunstan Dr , Toronto, Ontario M1L 2V5 Canada

\end

$$\matrix{T' & \longrightarrow & T \cr
  \downarrow q &              &\downarrow p \cr
B' &{\buildrel{F}\over{\longrightarrow}} & B}
$$

$B_nZ/_2$,  $B_n{\rm Z}/_2$,  $B_n{\rm Z}_{/2}$, $B_nZ/2$, $B_n{\rm Z}/2$